\documentclass[11pt,reqno]{amsart}
\usepackage[T1]{fontenc}
\usepackage{amsmath,amssymb,mathtools}
\usepackage{microtype,needspace}
\usepackage[hidelinks]{hyperref}
\newtheorem{maintheorem}{Theorem}

\newtheorem{lemma}{Lemma}
\theoremstyle{definition}
\newtheorem*{pietschproblem}{Problem 16}
\newcommand{\N}{\mathbb N}
\newcommand{\R}{\mathbb R}
\newcommand{\C}{\mathbb C}
\newcommand{\K}{\mathbb K}
\newcommand{\A}{\mathcal A}
\newcommand{\Lop}{\mathcal L}
\newcommand{\norm}[1]{\lVert #1\rVert}
\newcommand{\rank}{\operatorname{rank}}
\title[Pietsch's favourite problem on singular traces]{Pietsch's favourite problem on singular traces}
\author[G. Ara\'ujo]{Gustavo Ara\'ujo}
\address{Departamento de Matem\'atica, Universidade Estadual da Para\'iba, Campina Grande, PB, Brazil}
\email{gustavoaraujo@servidor.uepb.edu.br}
\author[D. Pellegrino]{Daniel Pellegrino}
\address{Departamento de Matem\'atica, Universidade Federal da Para\'iba, Jo\~ao Pessoa, PB, Brazil}
\email{daniel.pellegrino@academico.ufpb.br}
\author[A. Raposo Jr.]{Anselmo Raposo Jr.}
\address{Coordena\c{c}\~ao do Curso de Matem\'atica -- Bacharelado, Universidade Federal do Maranh\~ao, S\~ao Lu\'is, MA, Brazil}
\email{anselmo.junior@ufma.br}
\author[J. Santos]{Joedson Santos}
\address{Departamento de Matem\'atica, Universidade Federal da Para\'iba, Jo\~ao Pessoa, PB, Brazil}
\email{joedson.santos@academico.ufpb.br}
\begin{document}
\begin{abstract}
In 2023, A.~Pietsch formulated twenty-seven open problems concerning traces on operator ideals. We settle Problem~16, his personal favourite, by constructing nonzero bounded singular spectral traces on the weak Weyl ideal over the class of all Banach spaces. We also prove the asymptotic additivity of logarithmically normalized eigenvalue sums on this ideal.
\end{abstract}
\maketitle

\section{Introduction and main results}
In his 2023 paper on traces on operator ideals over the class of all Banach spaces, Pietsch \cite{Pietsch2023} asked whether the ideal of operators whose Weyl numbers decay at least as fast as the reciprocal of their index supports a nonzero trace (Problem~16). We settle this problem for both real and complex Banach spaces. We also show that the logarithmically normalized partial sums of the eigenvalues are additive up to a null sequence.
Let $\K=\R$ or $\C$, let $\N=\{1,2,\ldots\}$, $\N_0=\{0,1,\ldots\}$, and let $X,Y$ be Banach spaces over $\K$. Let $\Lop(X,Y)$ denote the Banach space of bounded linear operators from $X$ to $Y$. Write $X^*=\Lop(X,\K)$ and $I_X$ for the identity on $X$. Set
\begin{gather*}
 \rank T=\dim T(X)\quad(T\in\Lop(X,Y)),\\
 \mathcal F(X,Y)=\{T\in\Lop(X,Y):\rank T<\infty\}.
\end{gather*}
Operator norms are denoted by $\norm{\cdot}$.

For $F\in\mathcal F(X,X)$, write
\[
 F(x)=\sum_{j=1}^{m} f_j(x)x_j\quad(x\in X),
 \qquad f_j\in X^*,\quad x_j\in X.
\]
The trace $\operatorname{tr}\colon\mathcal F(X,X)\to\K$ is given by
\[
 \operatorname{tr}(F)=\sum_{j=1}^{m} f_j(x_j).
\]
It is independent of the representation of $F$, linear, and satisfies
$\operatorname{tr}(BA)=\operatorname{tr}(AB)$ for $A\in\mathcal F(X,Y)$ and $B\in\Lop(Y,X)$.

An operator ideal $\A$ over $\K$ assigns to each pair of Banach spaces $X,Y$ over $\K$ a linear subspace
\[
 \mathcal F(X,Y)\subseteq\A(X,Y)\subseteq\Lop(X,Y)
\]
such that, for all Banach spaces $E,G$ over $\K$,
\[
 V T U\in\A(E,G)\quad
 (T\in\A(X,Y),\ U\in\Lop(E,X),\ V\in\Lop(Y,G)).
\]
The finite-rank, compact, nuclear and absolutely summing operators are classical examples (see \cite{Pietsch1987}). A trace on $\A$ is a family of linear functionals $\tau_X\colon\A(X,X)\to\K$ satisfying
\begin{equation}\label{eq:trace}
 \tau_X(BA)=\tau_Y(AB)\quad(A\in\A(X,Y),\ B\in\Lop(Y,X)).
\end{equation}
For $X=Y$, \eqref{eq:trace} reduces to $\tau_X(BA-AB)=0$ for $A\in\A(X,X)$ and $B\in\Lop(X,X)$. The trace is called \emph{singular} if
\begin{equation}\label{eq:singular-def}
 \tau_X(F)=0\qquad(F\in\mathcal F(X,X))
\end{equation}
for every $X$. The theory of singular traces goes back to Dixmier \cite{Dixmier1966}. Pietsch \cite{Pietsch2012} developed an approach to Dixmier traces on the Marcinkiewicz operator ideal over Banach and Hilbert spaces using invariant functionals on sequence spaces. See also his historical account \cite{Pietsch2014History} and the monograph \cite{LordSukochevZanin2021} for the Hilbert space theory.

Write $\ell_2=\ell_2(\K)$, where
\[
 \ell_2(\K)=\left\{(u_j)_{j\ge1}\in\K^{\N}:\sum_{j=1}^\infty|u_j|^2<\infty\right\},
 \quad \norm u_2=\left(\sum_{j=1}^\infty|u_j|^2\right)^{1/2}.
\]
Let $(e_n)_{n\ge1}$ be the canonical basis of $\ell_2$. For $T\in\Lop(X,Y)$ and $n\in\N$, the $n$-th approximation number and the $n$-th Weyl number of $T$ are, respectively,
\begin{align}
 a_n(T)&=\inf\{\norm{T-F}:F\in\mathcal F(X,Y),\ \rank F<n\},\notag\\
 x_n(T)&=\sup\{a_n(TA):A\in\Lop(\ell_2,X),\ \norm A\le1\}.\label{eq:weyl-def}
\end{align}
In particular,
\begin{equation*}
 x_{n+1}(T)\le x_n(T)\le a_n(T)\qquad(n\in\N).
\end{equation*}
Taking $A=I_{\ell_2}$ in \eqref{eq:weyl-def}, we have
\begin{equation}\label{eq:hilbert-weyl}
 x_n(T)=a_n(T)\qquad(T\in\Lop(\ell_2,Y),\ n\in\N).
\end{equation}
For compact operators on $\ell_2$, $a_n(T)$ is the $n$-th singular value.

The weak Weyl ideal considered in \cite[Subsection~2.6]{Pietsch2023} is
\begin{align*}
 \norm T_{\mathcal W}&:=\sup_{n\ge1}n x_n(T),\\
 \mathcal W(X,Y)&:=\Lop_{1,\infty}^{\mathrm{WEYL}}(X,Y)
 =\{T\in\Lop(X,Y):\norm T_{\mathcal W}<\infty\}.
\end{align*}
It is a quasi-Banach operator ideal (see \cite{Pietsch2023}). For $S,T\in\mathcal W(X,Y)$ and $n\in\N$,
\begin{align}
 x_{2n-1}(S+T)&\le x_n(S)+x_n(T),\notag\\
 \norm{S+T}_{\mathcal W}&\le2(\norm S_{\mathcal W}+\norm T_{\mathcal W}).\label{eq:quasitr}
\end{align}
The dyadic quasi-norm in \cite{Pietsch2023} satisfies
\[
 \sup_{k\ge0}2^kx_{2^k}(T)\le\norm T_{\mathcal W}
 \le2\sup_{k\ge0}2^kx_{2^k}(T).
\]
A trace $\tau$ on $\mathcal W$ is bounded if there exists $C>0$, independent of $X$, such that
\[
 |\tau_X(T)|\le C\norm T_{\mathcal W}\qquad(T\in\mathcal W(X,X)).
\]

An operator between Banach spaces is Fredholm if its range is closed and both its kernel and the quotient of its codomain by its range are finite-dimensional. For a complex Banach space $X$, an operator $T\in\Lop(X,X)$ is a Riesz operator if, for every $\lambda\in\C\setminus\{0\}$, the range $R_\lambda=(\lambda I_X-T)(X)$ is closed and
\[
 \dim\ker(\lambda I_X-T)<\infty,\qquad \dim(X/R_\lambda)<\infty.
\]
By Lemma~\ref{lem:eigen}, every $T\in\mathcal W(X,X)$ is a Riesz operator. In particular, for $A\in\mathcal W(X,Y)$ and $B\in\Lop(Y,X)$,
\[
 BA\in\mathcal W(X,X),\qquad I_X-BA\text{ is Fredholm}.
\]
Thus $\mathcal W$ is an ideal of inessential operators in the sense of \cite[Subsection~2.4]{Pietsch2023}. The spectrum of an operator $T\in\Lop(X,X)$ over $\C$ is the set of $\lambda\in\C$ for which $\lambda I_X-T$ is not invertible in $\Lop(X,X)$. The nonzero spectral values of a Riesz operator are isolated eigenvalues of finite algebraic multiplicity, accumulating at most at zero. For $\lambda\ne0$, its algebraic multiplicity is the dimension of the generalized eigenspace $\bigcup_{j\ge1}\ker((T-\lambda I_X)^j)$ (zero if $\lambda$ is not an eigenvalue). We write $(\lambda_n(T))_{n\ge1}$ for the nonzero eigenvalues, repeated according to their algebraic multiplicities and ordered by nonincreasing modulus. If the list is finite, it is completed by zeros. When eigenvalues have the same modulus, choose any order among them; the normalized sums used in the main results will be independent of that choice modulo null sequences.
For the ideal $\mathcal W$ over $\C$, the definition of a \emph{spectral} trace in \cite[Subsection~2.4]{Pietsch2023} is
\begin{equation}\label{eq:spectral-def}
 \tau_X(T)=\varphi\left(\left(2^{-k}\sum_{n\in I_k}\lambda_n(T)\right)_{k\ge0}\right),
 \qquad I_k=\{2^k,\ldots,2^{k+1}-1\},
\end{equation}
where $\varphi\colon Z\to\C$ is a linear functional, independent of $X$, on the space
\begin{equation*}
 Z=\left\{b=(b_k)_{k\ge0}\in\C^{\N_0}:\sup_{k\ge0}2^k|b_k|<\infty\right\}.
\end{equation*}
For $b\in Z$, with $\ell_2=\ell_2(\C)$, define
\begin{equation}\label{eq:dyadicspace}
 B_b\colon\ell_2\to\ell_2,\qquad B_be_n=b_ke_n\quad(n\in I_k).
\end{equation}
Here $B_b\in\mathcal W(\ell_2,\ell_2)$: its diagonal entries $\beta_n=b_k$ for $n\in I_k$ satisfy $|\beta_n|\le 2M_b/n$, where $M_b=\sup_{k\ge0}2^k|b_k|$. The decreasing rearrangement of their moduli satisfies the same bound (at most $n-1$ entries can exceed $2M_b/n$). Thus its singular values, and hence its Weyl numbers, satisfy this bound. The traces constructed in this paper have $\varphi(b)=\tau_{\ell_2(\C)}(B_b)$. The sequence of dyadic means belongs to $Z$ by Lemma~\ref{lem:eigen}. For a real Banach space $X$, let $X_{\C}=X+iX$ be its complexification, endowed with the norm
\begin{equation}\label{eq:complex-norm}
 \norm{x+iy}_{X_{\C}}=\sup_{\theta\in\R}\norm{x\cos\theta-y\sin\theta}_X.
\end{equation}
For $T\in\Lop(X,Y)$ over $\R$, define
\begin{equation}\label{eq:complex-op}
 T_{\C}\colon X_{\C}\to Y_{\C},\qquad T_{\C}(x+iy)=Tx+iTy.
\end{equation}
By \eqref{eq:complexification}, $T_{\C}\in\mathcal W(X_{\C},X_{\C})$ for $T\in\mathcal W(X,X)$. Put $Z_{\R}=Z\cap\R^{\N_0}$. For $b\in Z_{\R}$, the same formula \eqref{eq:dyadicspace} defines a real diagonal operator on $\ell_2(\R)$. Over $\R$, a trace on $\mathcal W$ is spectral if there is a real linear functional $\varphi_{\R}\colon Z_{\R}\to\R$, independent of $X$, such that
\begin{equation*}
 \tau_X(T)=\varphi_{\R}\left(\left(2^{-k}\sum_{n\in I_k}\operatorname{Re}\lambda_n(T_{\C})\right)_{k\ge0}\right).
\end{equation*}

For Banach spaces $X,Y$ over $\C$, let
\begin{equation}\label{eq:strong-weyl}
 \mathcal V(X,Y)=\left\{T\in\Lop(X,Y):\sum_{n=1}^{\infty}x_n(T)<\infty\right\}.
\end{equation}
Since $nx_n(T)\le\sum_{j=1}^n x_j(T)$, $\mathcal V(X,Y)\subseteq\mathcal W(X,Y)$. By \cite[Corollary~2.3]{White1996}, the family
\begin{equation}\label{eq:strong-trace}
 \operatorname{tr}_{1,X}\colon\mathcal V(X,X)\to\C,
 \qquad \operatorname{tr}_{1,X}(T)=\sum_{n=1}^{\infty}\lambda_n(T),
\end{equation}
is a trace on $\mathcal V$. On the weak Weyl ideal $\mathcal W$, this sum need not converge. For instance, define
\begin{equation}\label{eq:diag}
 D\colon\ell_2\longrightarrow\ell_2,\qquad De_n=\frac1n e_n\quad(n\in\N).
\end{equation}
Then
\begin{equation}\label{eq:diagprops}
 x_n(D)=a_n(D)=\lambda_n(D)=\frac1n,
 \qquad \norm D_{\mathcal W}=1,
 \qquad \sum_{n=1}^{\infty}\lambda_n(D)=\infty.
\end{equation}
For $m\in\N$, define $D_m\colon\ell_2\to\ell_2$ by
\[
 D_me_n=\begin{cases}n^{-1}e_n,&1\le n\le m,\\0,&n>m.\end{cases}
\] Since
\[
 \rank D_m=m,\qquad \norm{D-D_m}=\frac1{m+1}\longrightarrow0,
\]
a singular trace taking a nonzero value at $D$ cannot be continuous in the operator norm.

By \eqref{eq:hilbert-weyl}, the weak trace class on the complex Hilbert space $\ell_2$ is
\[
 \Lop_{1,\infty}(\ell_2)=\mathcal W(\ell_2,\ell_2)
 =\{T\in\Lop(\ell_2,\ell_2):\sup_{n\ge1}n a_n(T)<\infty\}.
\] The spectral description of traces on Hilbert space ideals and the structure of traces on the weak trace class were studied in \cite{Semenov2015,SukochevZanin2014}. Let
\begin{gather*}
 \ell_\infty=\{u=(u_j)_{j\ge1}\in\C^{\N}:\norm u_\infty:=\sup_{j\ge1}|u_j|<\infty\},\\
 c_0=\{u\in\ell_\infty:\lim_{j\to\infty}u_j=0\}.
\end{gather*}
Write $[u]=u+c_0\in\ell_\infty/c_0$ for the class of $u\in\ell_\infty$, and let $\log$ denote the natural logarithm. Sequences indexed by $N\ge2$ are identified with sequences indexed by $\N$ by adjoining a zero first term. By \cite[Lemma~2.9]{Ponge2023}, the map
\begin{equation}\label{eq:hilbert-trace}
 \Lop_{1,\infty}(\ell_2)\longrightarrow\ell_\infty/c_0,\qquad
 T\longmapsto\left[\left(\frac{1}{\log N}\sum_{n=1}^N\lambda_n(T)\right)_{N\ge2}\right]
\end{equation}
is a bounded linear trace. An extended limit is a linear functional $\omega\colon\ell_\infty\to\C$ satisfying
\[
 \omega(u)\ge0\quad(u_j\ge0,\ j\in\N),\qquad
 \omega((1)_{j\ge1})=1,\qquad \omega|_{c_0}=0.
\]
By \cite[Proposition~2.10 and Definition~2.11]{Ponge2023}, the map
\begin{equation*}
 \operatorname{Tr}_\omega\colon\Lop_{1,\infty}(\ell_2)\to\C,\quad
 \operatorname{Tr}_\omega(T)=\omega\left(\left(\frac1{\log N}\sum_{n=1}^N\lambda_n(T)\right)_{N\ge2}\right)
\end{equation*}
is a Dixmier trace.

For ideals defined over the class of all Banach spaces, the cyclic identity \eqref{eq:trace} relates the components $\A(X,X)$ and $\A(Y,Y)$ for arbitrary $X$ and $Y$. Pietsch developed constructions of such traces using approximation numbers and generalized approximation numbers in \cite{Pietsch2015,Pietsch2017}. He also obtained nonzero traces on subideals of $\mathcal W$ defined by additional approximation conditions \cite[Lemma~2.5 and Theorem~2.8]{Pietsch2023}. For example, writing
\[
 \begin{gathered}
 \norm{T}_{\mathcal W_{2,\infty}}:=\sup_{n\ge1}\sqrt n\,x_n(T),\quad T\in\Lop(X,Y),\\
 \mathcal W_{2,\infty}(X,Y):=\{T\in\Lop(X,Y):\norm T_{\mathcal W_{2,\infty}}<\infty\},
 \end{gathered}
\]
there is a nonzero trace on the ideal of operators $T\in\Lop(X,Y)$ satisfying
\begin{equation*}
 \norm T_{\mathcal W_{2,\infty}}<\infty,\qquad
 \delta(T):=\sup_{n\ge1}\sqrt n\inf_{\substack{F\in\mathcal F(X,Y)\\\rank F<n}}\norm{T-F}_{\mathcal W_{2,\infty}}<\infty.
\end{equation*}
For $n\in\N$ and $F\in\mathcal F(X,Y)$ with $\rank F<n$,
\[
 x_{2n-1}(T)\le x_n(T-F)+x_n(F)=x_n(T-F)
 \le n^{-1/2}\norm{T-F}_{\mathcal W_{2,\infty}}.
\]
Taking the infimum over $F$,
\[
 n x_{2n-1}(T)\le\sqrt n\inf_{\substack{F\in\mathcal F(X,Y)\\\rank F<n}}\norm{T-F}_{\mathcal W_{2,\infty}}\le \delta(T).
\]
For $m\in\N$ and $n=\lfloor(m+1)/2\rfloor$,
\[
 m x_m(T)\le 2n x_{2n-1}(T)\le 2\delta(T).
\]
Thus $\norm T_{\mathcal W}\le 2\delta(T)<\infty$, which proves that the ideal in this example is contained in $\mathcal W$.

\Needspace{13\baselineskip}
The existence of a nonzero trace on $\mathcal W$ itself is Problem~16 in \cite[p.~30]{Pietsch2023}, introduced with the words ``The following problem is my personal favourite.''
\begin{pietschproblem}[{\cite[p.~30]{Pietsch2023}}]
Show that $\Lop_{1,\infty}^{\mathrm{WEYL}}$ supports non-vanishing traces.
\end{pietschproblem}

\begin{maintheorem}\label{thm:weyltrace}
For each scalar field $\K=\R$ or $\C$, the weak Weyl ideal $\mathcal W$ over the class of all Banach spaces over $\K$ admits a bounded singular spectral trace $\tau=(\tau_X)_X$ such that $\tau_{\ell_2(\K)}(D)=1$, where $D$ is the diagonal operator in \eqref{eq:diag} over $\K$.
\end{maintheorem}

The Hilbert space result \eqref{eq:hilbert-trace} also extends to arbitrary Banach spaces. For complex $X$, put $\mu_n(T)=\lambda_n(T)$; for real $X$, put $\mu_n(T)=\operatorname{Re}\lambda_n(T_{\C})$, with the eigenvalues in both cases ordered by nonincreasing modulus before taking real parts. Write $\ell_\infty(\K)$ and $c_0(\K)$ for the bounded and null sequences over $\K$, respectively, with the same conventions as in the complex case. The quotient carries the norm $\norm{[u]}=\inf_{v\in c_0(\K)}\norm{u-v}_\infty$.

\begin{maintheorem}\label{thm:quotient}
For every Banach space $X$ over $\K=\R$ or $\C$, the map
\begin{equation}\label{eq:quotient-trace}
 \begin{aligned}
 \Psi_X&\colon\mathcal W(X,X)\longrightarrow\ell_\infty(\K)/c_0(\K),\\
 \Psi_X(T)&=\left[\left(\frac1{\log N}\sum_{n=1}^N\mu_n(T)\right)_{N\ge2}\right].
 \end{aligned}
\end{equation}
is well defined, independently of the ordering at ties in modulus, and linear. It satisfies
\begin{align*}
 \Psi_X(BA)&=\Psi_Y(AB)\quad(A\in\mathcal W(X,Y),\ B\in\Lop(Y,X)),\\
 \Psi_X(F)&=0\quad(F\in\mathcal F(X,X)),\qquad \Psi_{\ell_2}(D)=[\mathbf 1],
\end{align*}
where $\mathbf 1$ is the constant sequence with value $1$. Moreover,
\[
 \norm{\Psi_X(T)}_{\ell_\infty(\K)/c_0(\K)}\le C_{\K}\norm T_{\mathcal W},
 \qquad C_{\C}=2\sqrt2\,\mathrm e^2,\quad C_{\R}=8\sqrt2\,\mathrm e^2.
\]
\end{maintheorem}

In particular, for every $S,T\in\mathcal W(X,X)$,
\begin{equation}\label{eq:asymptotic-additivity}
 \lim_{N\to\infty}\frac{1}{\log N}
 \sum_{n=1}^N\big(\mu_n(S+T)-\mu_n(S)-\mu_n(T)\big)=0.
\end{equation}
By Theorem~\ref{thm:quotient}, every extended limit (restricted to real sequences in the real case) yields a bounded singular spectral trace on $\mathcal W$, normalized by its value at $D$.

The main difficulty in the proof is the additivity of the logarithmically normalized eigenvalue sums. Although White's trace \eqref{eq:strong-trace} is not defined on all of $\mathcal W$, products of two operators in $\mathcal W$ belong to $\mathcal V$. For $S,T\in\mathcal W(X,X)$ over $\C$, the spectral traces of powers of $S+zT$, $z\in\C$, determine a regularized determinant holomorphic in $z$. A Fourier coefficient of the logarithm of its modulus approximates the normalized eigenvalue sums with an error of order $1/\log N$ for $N\ge2$, locally uniformly in $z$. Jensen's formula and the submean inequality give a mean-value defect of the same order. The ultralimits are therefore harmonic; their linear growth and complex homogeneity yield additivity. Finally, the arbitrariness of the ultrafilter gives \eqref{eq:asymptotic-additivity}.

\section{Auxiliary results}
Throughout this section $X$ is a complex Banach space.

\begin{lemma}\label{lem:eigen}
Every $T\in\mathcal W(X,X)$ is a Riesz operator and
\begin{equation}\label{eq:weak-eig}
 |\lambda_m(T)|\le\frac{C_0\norm T_{\mathcal W}}{m}\quad(m\in\N),
 \qquad C_0=2\sqrt2\,\mathrm e^2.
\end{equation}
\end{lemma}
\begin{proof}
Suppose that there exist an infinite-dimensional closed subspace $E\subset X$ and $c>0$ such that
\begin{equation}\label{eq:lower-bd}
 \norm{Tx}_X\ge c\norm x_X\qquad(x\in E).
\end{equation} By Dvoretzky's theorem \cite{Dvoretzky1961}, for each $n\in\N$, writing $\ell_2^n=(\C^n,\norm{\cdot}_2)$, there exists $J_n\colon\ell_2^n\to E$ such that
\[
 \tfrac12\norm u_{\ell_2^n}\le\norm{J_nu}_X\le\norm u_{\ell_2^n}\quad(u\in\ell_2^n).
\]
Let $P_n\colon\ell_2\to\ell_2^n$, $P_n((u_j)_{j\ge1})=(u_1,\ldots,u_n)$. For $F\in\mathcal F(\ell_2,X)$ with $\rank F<n$, choose $u\in\ker F\cap\operatorname{span}\{e_1,\ldots,e_n\}$ with $\norm u_2=1$. By \eqref{eq:lower-bd},
\[
 \norm{TJ_nP_n-F}\ge\norm{TJ_nP_nu}_X
 \ge c\norm{J_nP_nu}_X\ge c/2.
\]
Taking the infimum over $F$ and using $\norm{J_nP_n}\le1$, we obtain
\[
 \frac{\norm T_{\mathcal W}}{n}\ge x_n(T)\ge a_n(TJ_nP_n)\ge c/2
 \qquad(n\in\N),
\]
a contradiction. Therefore $T$ is strictly singular (it is not bounded below on any infinite-dimensional closed subspace of $X$) and, by \cite[p.~66]{Caradus1966}, a Riesz operator.

By the Weyl eigenvalue inequality \cite{Pietsch1980Weyl},
\begin{align*}
 |\lambda_{2n-1}(T)|
 &\le\sqrt2\,\mathrm e\left(\prod_{j=1}^n x_j(T)\right)^{1/n}\\
 &\le\sqrt2\,\mathrm e\norm T_{\mathcal W}(n!)^{-1/n}
 \le\frac{\sqrt2\,\mathrm e^2\norm T_{\mathcal W}}{n},
\end{align*}
where we used
\[
 \log(n!)\ge\int_1^n\log t\,dt=n\log n-n+1.
\]
For $m\in\N$, take $n=\lfloor(m+1)/2\rfloor$. Then
\[
 m|\lambda_m(T)|\le 2n|\lambda_{2n-1}(T)|\le2\sqrt2\,\mathrm e^2\norm T_{\mathcal W}.
\]
\end{proof}

For $T\in\mathcal W(X,X)$ and two orderings $(\lambda_n(T))_{n\ge1}$ and $(\lambda'_n(T))_{n\ge1}$ differing only at ties in modulus, the first $N$ terms can differ as multisets only in eigenvalues of modulus $|\lambda_N(T)|$. By \eqref{eq:weak-eig},
\begin{equation}\label{eq:ties}
 \left|\sum_{n=1}^N(\lambda_n(T)-\lambda'_n(T))\right|
 \le2N|\lambda_N(T)|\le2C_0\norm T_{\mathcal W}.
\end{equation}

A free ultrafilter on $\N$ is a family of subsets of $\N$ closed under finite intersections and supersets, not containing the empty set, containing every cofinite set, and containing exactly one of $E$ and $\N\setminus E$ for every $E\subseteq\N$. Fix an arbitrary free ultrafilter $\mathcal U$ on $\N$. For every $u=(u_N)_{N\ge1}\in\ell_\infty$, there is a unique $b\in\C$, denoted by $\lim_{N\to\mathcal U}u_N$, such that
\[
 \{N\in\N:|u_N-b|<\varepsilon\}\in\mathcal U\qquad(\varepsilon>0).
\]
The map $\omega_{\mathcal U}\colon\ell_\infty\to\C$, $\omega_{\mathcal U}(u)=\lim_{N\to\mathcal U}u_N$, is an extended limit: ultralimits preserve linear combinations and positivity, and agree with ordinary limits when the latter exist. For $N\ge2$,
\begin{equation}\label{eq:harmonic-sum}
 \log(N+1)\le\sum_{n=1}^N\frac1n\le1+\log N.
\end{equation}
By \eqref{eq:harmonic-sum} and \eqref{eq:weak-eig},
\[
 \left|\frac1{\log N}\sum_{n=1}^N\lambda_n(T)\right|
 \le C_0\norm T_{\mathcal W}\frac{1+\log N}{\log N}\qquad(N\ge2).
\]
For each complex Banach space $X$, define the candidate trace $\tau_X\colon\mathcal W(X,X)\to\C$ by
\begin{equation}\label{eq:singulartrace}
 \tau_X(T)=\lim_{N\to\mathcal U}\frac1{\log N}\sum_{n=1}^N\lambda_n(T),
 \qquad T\in\mathcal W(X,X),
\end{equation}
where $N\ge2$. By \eqref{eq:ties}, the value is independent of the ordering at ties in modulus. On the complex space $\ell_2$, \eqref{eq:singulartrace} agrees with the Dixmier trace in \cite[Proposition~2.10 and Definition~2.11]{Ponge2023} associated with the ultralimit.

For $T\in\mathcal W(X,X)$, define the regularized determinant by
\begin{equation}\label{eq:det2}
 d_T(w)=\prod_{n=1}^\infty (1-w\lambda_n(T))\exp(w\lambda_n(T)),\qquad w\in\C.
\end{equation}
Since $\sum_{n=1}^\infty|\lambda_n(T)|^2<\infty$, the product converges locally uniformly and defines an entire function $d_T\colon\C\to\C$.

\begin{lemma}\label{lem:detanalytic}
If $S,T\in\mathcal W(X,X)$, then $(z,w)\mapsto d_{S+zT}(w)$ is holomorphic on $\C^2$.
\end{lemma}
\begin{proof}

By multiplicativity and monotonicity of Weyl numbers \cite{Pietsch1980Weyl},
\begin{align*}
 x_{2n}(AB)&\le x_{2n-1}(AB)\le x_n(A)x_n(B)\\
 &\le n^{-2}\norm A_{\mathcal W}\norm B_{\mathcal W}\qquad(A,B\in\mathcal W(X,X)).
\end{align*}
Consequently,
\begin{equation}\label{eq:product}
 \sum_{j=1}^{\infty}x_j(AB)\le2\norm A_{\mathcal W}\norm B_{\mathcal W}\sum_{n=1}^{\infty}\frac1{n^2}<\infty.
\end{equation}
By \eqref{eq:strong-weyl}, \eqref{eq:product} and the ideal property of $\mathcal V$,
\[
 U_1\cdots U_m\in\mathcal V(X,X)
 \qquad(U_1,\ldots,U_m\in\mathcal W(X,X),\ m\ge2).
\]
Fix $m\ge2$ and $z\in\C$ and put $Q=S+zT$. For $\nu\in\C\setminus\{0\}$ and $k\in\N$, the polynomials $(t-\lambda)^k$, with $\lambda^m=\nu$, are pairwise relatively prime. By the kernel decomposition for relatively prime polynomials,
\[
 \ker((Q^m-\nu I_X)^k)
 =\bigoplus_{\lambda^m=\nu}\ker((Q-\lambda I_X)^k).
\]
Since $Q$ and $Q^m$ are Riesz operators, for sufficiently large $k$ these kernels are the corresponding generalized eigenspaces (or $\{0\}$). For a Riesz operator $R\in\Lop(X,X)$ and $\lambda\in\C\setminus\{0\}$, write
\[
 \operatorname{mult}(\lambda,R)=\dim\left(\bigcup_{j\ge1}\ker((R-\lambda I_X)^j)\right).
\]
Then
\[
 \operatorname{mult}(\nu,Q^m)=\sum_{\lambda^m=\nu}\operatorname{mult}(\lambda,Q).
\]
Thus, by \eqref{eq:strong-trace},
\begin{equation}\label{eq:powertrace}
 p_m(z)=\sum_{n=1}^{\infty}\lambda_n(S+zT)^m
 =\operatorname{tr}_{1,X}((S+zT)^m).
\end{equation}
Let $U_0=S$, $U_1=T$. By linearity of $\operatorname{tr}_{1,X}$,
\[
 p_m(z)=\sum_{\varepsilon\in\{0,1\}^m}
 z^{\varepsilon_1+\cdots+\varepsilon_m}
 \operatorname{tr}_{1,X}(U_{\varepsilon_1}\cdots U_{\varepsilon_m}).
\]
Thus $p_m\colon\C\to\C$ is a polynomial of degree at most $m$.

Fix a closed disk $K\subset\C$. Put $M=\max\{1,2C_0(\norm S_{\mathcal W}+\sup_{z\in K}|z|\norm T_{\mathcal W})\}$. By \eqref{eq:weak-eig} and \eqref{eq:quasitr},
\begin{equation}\label{eq:uniformeig}
 |\lambda_n(S+zT)|\le M/n\quad(z\in K,\ n\in\N).
\end{equation}
For $z\in K$ and $|w|<1/M$, \eqref{eq:det2} and \eqref{eq:powertrace} yield
\begin{equation}\label{eq:detseries}
 d_{S+zT}(w)=\exp\left(-\sum_{m=2}^{\infty}\frac{w^mp_m(z)}m\right).
\end{equation}
For $0<r<1/M$, by \eqref{eq:powertrace} and \eqref{eq:uniformeig},
\begin{equation}\label{eq:normal-series}
 \sum_{m=2}^{\infty}\sup_{\substack{z\in K\\|w|\le r}}\frac{|w^mp_m(z)|}{m}
 \le\sum_{m=2}^{\infty}\frac{(rM)^m}{m}\sum_{n=1}^{\infty}n^{-2}<\infty.
\end{equation}
By \eqref{eq:detseries} and \eqref{eq:normal-series}, the Taylor coefficients
\[
 c_j(z)=\frac1{j!}\left.\frac{\partial^j}{\partial w^j}d_{S+zT}(w)\right|_{w=0}\quad(j\in\N_0)
\]
are entire functions of $z$.

For $|u|\le1/2$,
\begin{equation}\label{eq:factor-small}
 \big|\log|(1-u)e^u|\big|\le\sum_{m=2}^{\infty}\frac{|u|^m}{m}
 \le\frac{|u|^2}{2(1-|u|)}\le |u|^2.
\end{equation}
For $|u|>1/2$,
\begin{equation}\label{eq:factor-large}
 \log|(1-u)e^u|\le\log(1+|u|)+|u|\le2|u|\le4|u|^2.
\end{equation}
Applying the bounds \eqref{eq:factor-small} and \eqref{eq:factor-large} to each factor of the determinant and using \eqref{eq:uniformeig}, we obtain
\begin{equation}\label{eq:detbound}
 |d_{S+zT}(w)|\le\exp\left(4|w|^2M^2\sum_{n=1}^{\infty}n^{-2}\right).
\end{equation}
For any $R>0$, Cauchy's estimate and \eqref{eq:detbound} yield
\begin{equation}\label{eq:cauchy-coef}
 \sup_{z\in K}|c_j(z)|\le R^{-j}\exp\left(4R^2M^2\sum_{n=1}^{\infty}n^{-2}\right).
\end{equation}
For $|w|\le R/2$, by \eqref{eq:cauchy-coef},
\[
 \sum_{j=0}^{\infty}\sup_{z\in K}|c_j(z)w^j|
 \le 2\exp\left(4R^2M^2\sum_{n=1}^{\infty}n^{-2}\right).
\]
Consequently, $(z,w)\mapsto\sum_{j=0}^{\infty}c_j(z)w^j=d_{S+zT}(w)$ is holomorphic on $\C^2$.

\end{proof}

\begin{lemma}\label{lem:harmonic}
Let $S,T\in\mathcal W(X,X)$ and let $\tau_X$ be defined by \eqref{eq:singulartrace}. The function $z\mapsto\tau_X(S+zT)$ is harmonic on $\C$ (in its real and imaginary parts).
\end{lemma}
\begin{proof}
Let $P\colon\C\to\mathcal W(X,X)$ be given by $P(z)=S+zT$. For $N\ge2$, define $m_N\colon\C\to\R$ and $g_N\colon\C\to\C$ by
\begin{align}
 m_N(z)&=\frac1{2\pi}\int_0^{2\pi}\log|d_{P(z)}(Ne^{i\theta})|\,d\theta,\label{eq:jensenmean}\\
 g_N(z)&=\frac{1}{\pi N\log N}\int_0^{2\pi}e^{-i\theta}\log|d_{P(z)}(Ne^{i\theta})|\,d\theta.\label{eq:fourier}
\end{align}

For fixed $z$, $d_{P(z)}(0)=1$, so $d_{P(z)}$ has only finitely many zeros on the circle $|w|=N$. The logarithmic singularities at these zeros are integrable. Set $\log^+t=\max\{0,\log t\}$ for $t>0$ and $\log^+0=0$. Jensen's formula yields
\begin{equation}\label{eq:jensen}
 m_N(z)=\sum_{n=1}^{\infty}\log^+(N|\lambda_n(P(z))|).
\end{equation}
For $a\in\C$ with $|a|<1$,
\[
 \log|1-ae^{i\theta}|=-\operatorname{Re}\sum_{j=1}^{\infty}\frac{a^je^{ij\theta}}j.
\]
If $|a|>1$, then
\[
 \log|1-ae^{i\theta}|=\log|a|-\operatorname{Re}\sum_{j=1}^{\infty}\frac{a^{-j}e^{-ij\theta}}j.
\]
Integrating against $e^{-i\theta}$, and using
\[
 \frac1{2\pi}\int_0^{2\pi}e^{-i\theta}\operatorname{Re}(ae^{i\theta})\,d\theta=\frac a2,
\]
we obtain
\begin{equation}\label{eq:factorfourier}
 \frac1{2\pi}\int_0^{2\pi}e^{-i\theta}\log|(1-ae^{i\theta})e^{ae^{i\theta}}|\,d\theta
 =\begin{cases}0,& |a|\le1,\\ (a-1/\bar a)/2,& |a|>1.\end{cases}
\end{equation}
The boundary case follows by $L_1$-convergence of the logarithms as $ra\to a$, $r\uparrow1$: after a rotation, for $r\ge1/2$ and $-\pi\le\theta\le\pi$,
\[
 |1-re^{i\theta}|^2=(1-r)^2+4r\sin^2(\theta/2)\ge2\theta^2/\pi^2,
\]
and $|1-re^{i\theta}|\le2$. Dominated convergence applies since $|\log|\theta||$ is integrable on this interval. By \eqref{eq:factor-small}, for $N|\lambda_n(P(z))|\le1/2$,
\[
 \sup_{\theta}\left|\log|(1-Ne^{i\theta}\lambda_n(P(z)))e^{Ne^{i\theta}\lambda_n(P(z))}|\right|
 \le N^2|\lambda_n(P(z))|^2.
\]
By \eqref{eq:weak-eig} the series of these bounds converges. Applying \eqref{eq:factorfourier} to \eqref{eq:det2} and \eqref{eq:fourier} yields
\begin{equation}\label{eq:cutoff}
 g_N(z)=\frac1{\log N}\sum_{n:\,|\lambda_n(P(z))|>1/N}
 \left(\lambda_n(P(z))-\frac1{N^2\overline{\lambda_n(P(z))}}\right).
\end{equation}

Fix a closed disk $K\subset\C$ and set
\[
 M=\max\{1,2C_0(\norm S_{\mathcal W}+\sup_{z\in K}|z|\norm T_{\mathcal W})\}.
\] By \eqref{eq:weak-eig} and \eqref{eq:quasitr}, the bound \eqref{eq:uniformeig} holds for this $K$ and $M$. For $z\in K$ and $t>0$, put
\begin{equation*}
 q_z(t)=\#\{n\in\N:|\lambda_n(P(z))|>t\}\le M/t,
\end{equation*}
where $\#$ denotes cardinality. The inequality follows from \eqref{eq:uniformeig}.
Writing $\mathbf1_E$ for the indicator function of a set $E$, we obtain by \eqref{eq:jensen} and Tonelli's theorem
\begin{align}\label{eq:meanbound}
 0\le m_N(z)&=\sum_{n=1}^{\infty}\int_{1/N}^{\infty}
 \mathbf1_{\{t<|\lambda_n(P(z))|\}}\,\frac{dt}{t}\notag\\
 &=\int_{1/N}^{\infty}q_z(t)\,\frac{dt}{t}
 \le M\int_{1/N}^{\infty}\frac{dt}{t^2}=MN.
\end{align}
For $N\ge2$, define
\begin{equation}\label{eq:partial-mean}
 f_N\colon\C\to\C,\qquad f_N(z)=\frac1{\log N}\sum_{n=1}^N\lambda_n(P(z)).
\end{equation}
Fix $z\in K$ and write $q=q_z(1/N)$.
If $q\le N$, then
\begin{equation}\label{eq:short-cut}
 \left|\sum_{n=q+1}^N\lambda_n(P(z))\right|\le\frac{N-q}N\le1.
\end{equation}
If $q>N$, then $q\le MN$ and
\begin{equation}\label{eq:long-cut}
 \left|\sum_{n=N+1}^q\lambda_n(P(z))\right|
 \le M\sum_{n=N+1}^q\frac1n\le M\int_N^{MN}\frac{dt}{t}=M\log M.
\end{equation}
In either case,
\begin{equation}\label{eq:reciprocal-cut}
 \sum_{n=1}^q\frac1{N^2|\lambda_n(P(z))|}\le q/N\le M.
\end{equation}
By \eqref{eq:short-cut}, \eqref{eq:long-cut} and \eqref{eq:reciprocal-cut},
\begin{equation}\label{eq:approxsum}
 \sup_{z\in K}|f_N(z)-g_N(z)|\le\frac{1+M+M\log M}{\log N}.
\end{equation}
By \eqref{eq:uniformeig},
\[
 \sup_{z\in K}|f_N(z)|\le\frac{M}{\log N}\sum_{n=1}^N\frac1n
 \le M\frac{1+\log N}{\log N}.
\]
By \eqref{eq:approxsum},
\begin{equation}\label{eq:uniform-g-bound}
 L:=\sup_{N\ge2}\sup_{z\in K}|g_N(z)|<\infty.
\end{equation}

Let $\overline{B}(z,r)=\{\zeta\in\C:|\zeta-z|\le r\}\subset K$, with $r>0$. For each $\theta$, Lemma~\ref{lem:detanalytic} implies that $\log|d_{P(\cdot)}(Ne^{i\theta})|$ is subharmonic (or identically $-\infty$). For fixed $z$, the exceptional angles for which $d_{P(z)}(Ne^{i\theta})=0$ form a finite set, since $d_{P(z)}(0)=1$. For all other angles the subharmonic function is not identically $-\infty$. The submean inequality yields, for almost every $\theta\in[0,2\pi]$,
\begin{equation}\label{eq:submean-positive}
 \Delta_N(\theta):=\frac1{2\pi}\int_0^{2\pi}\log|d_{P(z+re^{it})}(Ne^{i\theta})|\,dt-\log|d_{P(z)}(Ne^{i\theta})|\ge0.
\end{equation}

Put $u_N(\zeta,\theta)=\log|d_{P(\zeta)}(Ne^{i\theta})|$ for $\zeta\in K$ and $\theta\in[0,2\pi]$. Write $u_N^+=\max\{u_N,0\}$ and $u_N^-=\max\{-u_N,0\}$ for its positive and negative parts. By \eqref{eq:detbound}, $u_N^+\le H_N$, where $H_N=4N^2M^2\sum_{n=1}^{\infty}n^{-2}$. By \eqref{eq:jensenmean} and \eqref{eq:meanbound},
\[
 \int_0^{2\pi}u_N^-(\zeta,\theta)\,d\theta
 =\int_0^{2\pi}u_N^+(\zeta,\theta)\,d\theta-2\pi m_N(\zeta)
 \le2\pi H_N.
\]
In particular,
\[
 \int_0^{2\pi}\int_0^{2\pi}|u_N(z+re^{it},\theta)|\,d\theta\,dt\le8\pi^2H_N<\infty.
\]
By Fubini's theorem, \eqref{eq:submean-positive}, \eqref{eq:fourier} and \eqref{eq:meanbound},

\begin{align}
 &\left|\frac1{2\pi}\int_0^{2\pi}g_N(z+re^{it})\,dt-g_N(z)\right|\notag\\
 &\qquad\le\frac{2}{N\log N}\frac1{2\pi}\int_0^{2\pi}\Delta_N(\theta)\,d\theta\notag\\
 &\qquad=\frac{2}{N\log N}\left(\frac1{2\pi}\int_0^{2\pi}m_N(z+re^{it})\,dt-m_N(z)\right)\notag\\
 &\qquad\le\frac{2M}{\log N}.\label{eq:mean-defect}
\end{align}

By Lemma~\ref{lem:detanalytic}, $(\zeta,\theta)\mapsto d_{P(\zeta)}(Ne^{i\theta})$ is continuous. With $\log 0=-\infty$, the logarithm of its modulus is jointly measurable. Set the integrand in \eqref{eq:fourier} equal to zero on the zero set of $d_{P(\zeta)}(Ne^{i\theta})$. It is jointly measurable and integrable in $\theta$ for each $\zeta\in\C$; thus $g_N$ is measurable. Let $K_0\subset\operatorname{int} K$ be a compact disk and choose a nonnegative, smooth, radially symmetric function $\rho\colon\C\to\R$ of compact support such that $\int_{\C}\rho\,d\mu=1$, where $\mu$ is planar Lebesgue measure on $\C\simeq\R^2$, $\nabla$ is the real gradient, and $\norm{\nabla\rho}_{L_1(\C)}=\int_{\C}|\nabla\rho(w)|\,d\mu(w)$. Take $\operatorname{supp}\rho\subset\{w:|w|<\operatorname{dist}(K_0,\partial K)\}$. Extend $g_N|_K$ by zero outside $K$, denoting the extension by $\widehat g_N$. Define
\[
 h_N\colon\C\to\C,\qquad h_N(z)=\int_{\C}\widehat g_N(w)\rho(z-w)\,d\mu(w).
\]
Write $\rho(w)=\rho_0(|w|)$, where $\rho_0\colon[0,\infty)\to[0,\infty)$. For $z\in K_0$, the change of variables $w=z+re^{it}$ gives
\[
 h_N(z)=\int_0^\infty\rho_0(r)r\int_0^{2\pi}g_N(z+re^{it})\,dt\,dr.
\]
By \eqref{eq:mean-defect},
\begin{align}
 |h_N(z)-g_N(z)|
 &\le2\pi\int_0^{\infty}\rho_0(r)r
 \left|\frac1{2\pi}\int_0^{2\pi}g_N(z+re^{it})\,dt-g_N(z)\right|dr\notag\\
 &\le\frac{2M}{\log N}\,2\pi\int_0^{\infty}\rho_0(r)r\,dr
 =\frac{2M}{\log N}.\label{eq:mollify}
\end{align}
By \eqref{eq:uniform-g-bound}, for $z,z'\in K_0$,
\begin{align}
 |h_N(z)-h_N(z')|
 &\le L\int_{\C}|\rho(z-w)-\rho(z'-w)|\,d\mu(w)\notag\\
 &\le L|z-z'|\int_0^1\int_{\C}|\nabla\rho(z'+t(z-z')-w)|\,d\mu(w)\,dt\notag\\
 &= L\norm{\nabla\rho}_{L_1(\C)}|z-z'|.\label{eq:mollipschitz}
\end{align}
Define $g\colon\C\to\C$ by $g(z)=\lim_{N\to\mathcal U}g_N(z)$. By \eqref{eq:approxsum} and \eqref{eq:singulartrace},
\begin{equation}\label{eq:glimit}
 g(z)=\tau_X(P(z))\qquad(z\in\C).
\end{equation}
By \eqref{eq:mollify} and \eqref{eq:mollipschitz},
\begin{equation}\label{eq:equicont}
 |g_N(z)-g_N(z')|\le\frac{4M}{\log N}+L\norm{\nabla\rho}_{L_1(\C)}|z-z'|.
\end{equation}
Taking ultralimits in \eqref{eq:equicont}, we have
\begin{equation}\label{eq:g-lip}
 |g(z)-g(z')|\le L\norm{\nabla\rho}_{L_1(\C)}|z-z'|\qquad(z,z'\in K_0).
\end{equation} For $\delta>0$, choose a finite $\delta$-net $\{z_1,\ldots,z_l\}\subset K_0$ (every point of $K_0$ is within distance $\delta$ of at least one $z_j$). Then
\begin{align*}
 \sup_{z\in K_0}|g_N(z)-g(z)|
 &\le\max_{1\le j\le l}|g_N(z_j)-g(z_j)|\\
 &\quad+\frac{4M}{\log N}+2L\norm{\nabla\rho}_{L_1(\C)}\delta.
\end{align*}
The ultralimit of the maximum is zero. Letting $\delta\downarrow0$, we conclude that
\begin{equation}\label{eq:uniformultra}
 \lim_{N\to\mathcal U}\sup_{z\in K_0}|g_N(z)-g(z)|=0.
\end{equation}
Since $K_0$ is arbitrary, \eqref{eq:mean-defect} and \eqref{eq:uniformultra} imply
\[
 g(z)=\frac1{2\pi}\int_0^{2\pi}g(z+re^{it})\,dt\qquad(z\in\C,\ r>0).
\]
By \eqref{eq:g-lip}, the mean-value characterization and \eqref{eq:glimit}, $z\mapsto\tau_X(S+zT)$ is harmonic.

\end{proof}

\section{Proof of Theorem A}
\begin{proof}

Let $X$ be a complex Banach space and $T\in\mathcal W(X,X)$. By \eqref{eq:weak-eig} and \eqref{eq:harmonic-sum},
\begin{equation}\label{eq:tracebound}
 |\tau_X(T)|\le C_0\norm T_{\mathcal W}\qquad(T\in\mathcal W(X,X)).
\end{equation}
Since $(\alpha\lambda_n(T))_{n\ge1}$ is an eigenvalue sequence of $\alpha T$,
\begin{equation}\label{eq:homog}
 \tau_X(\alpha T)=\alpha\tau_X(T)\quad(\alpha\in\C,\ T\in\mathcal W(X,X)).
\end{equation}

For $S,T\in\mathcal W(X,X)$, let $f\colon\C\to\C$ be given by $f(z)=\tau_X(S+zT)$. By Lemma~\ref{lem:harmonic}, $f$ is harmonic. By \eqref{eq:tracebound} and \eqref{eq:quasitr},
\begin{equation}\label{eq:linear-growth}
 |f(z)|\le2C_0(\norm S_{\mathcal W}+|z|\norm T_{\mathcal W})
 \le C_{S,T}(1+|z|)\quad(z\in\C),
\end{equation}
where $C_{S,T}=2C_0(\norm S_{\mathcal W}+\norm T_{\mathcal W})$.
By \eqref{eq:linear-growth} and the interior second derivative estimate for harmonic functions, at each $z_0\in\C$,
\[
 \norm{\nabla^2 f(z_0)}\le\frac{C}{R^2}\sup_{|z-z_0|\le R}|f(z)|
 \le\frac{CC_{S,T}(1+|z_0|+R)}{R^2}\xrightarrow[R\to\infty]{}0,
\]
where $C$ is an absolute constant and $\nabla^2 f$ is the real Hessian (componentwise). Thus $\nabla^2 f=0$, and there exist $\eta,\xi\in\C$ such that
\begin{equation}\label{eq:affine}
 f(z)=\tau_X(S)+\eta z+\xi\bar z\qquad(z\in\C).
\end{equation} By homogeneity \eqref{eq:homog} and the continuity of the harmonic function $w\mapsto\tau_X(T+wS)$ (Lemma~\ref{lem:harmonic}),
\begin{equation}\label{eq:asymptotic-slope}
 \frac{f(z)}z=\tau_X(T+z^{-1}S)\xrightarrow[|z|\to\infty]{}\tau_X(T).
\end{equation}
Taking $z=r$ and $z=ir$ in \eqref{eq:affine} and \eqref{eq:asymptotic-slope}, we obtain
\[
 \eta+\xi=\lim_{r\to\infty}\frac{f(r)}r=\tau_X(T)
 =\lim_{r\to\infty}\frac{f(ir)}{ir}=\eta-\xi.
\]
Thus $\xi=0$, $\eta=\tau_X(T)$ and
\begin{equation}\label{eq:additivity}
 \tau_X(S+T)=f(1)=\tau_X(S)+\tau_X(T).
\end{equation}

Let $Y$ be a complex Banach space, $A\in\mathcal W(X,Y)$, and $B\in\Lop(Y,X)$. For $\lambda\in\C\setminus\{0\}$ and $m\in\N$,
\begin{align*}
 A(BA-\lambda I_X)^m&=(AB-\lambda I_Y)^m A,\\
 B(AB-\lambda I_Y)^m&=(BA-\lambda I_X)^m B.
\end{align*}
For sufficiently large $m$, set
\begin{align*}
 E_\lambda&=\ker((BA-\lambda I_X)^m), & G_\lambda&=\ker((AB-\lambda I_Y)^m).
\end{align*}
Then $A(E_\lambda)\subseteq G_\lambda$ and $B(G_\lambda)\subseteq E_\lambda$. Since $BA|_{E_\lambda}$ and $AB|_{G_\lambda}$ are invertible, the restrictions $A|_{E_\lambda}$ and $B|_{G_\lambda}$ are injective. Hence
\[
 \dim E_\lambda\le\dim G_\lambda\le\dim E_\lambda.
\] Thus $BA$ and $AB$ have the same nonzero eigenvalues with the same algebraic multiplicities. By \eqref{eq:ties} and \eqref{eq:singulartrace}, $\tau_X(BA)=\tau_Y(AB)$. If $F\in\mathcal F(X,X)$, then $\sum_{n=1}^\infty|\lambda_n(F)|<\infty$ and
\[
 \left|\frac1{\log N}\sum_{n=1}^N\lambda_n(F)\right|
 \le\frac1{\log N}\sum_{n=1}^{\infty}|\lambda_n(F)|\longrightarrow0.
\]
Thus \eqref{eq:singular-def} holds. By \eqref{eq:diagprops},
\begin{equation}\label{eq:normalization}
 \tau_{\ell_2}(D)=\lim_{N\to\infty}\frac1{\log N}\sum_{n=1}^N\frac1n=1.
\end{equation}

For $T\in\mathcal W(X,X)$ and $k\in\N_0$, let
\[
 b_k=2^{-k}\sum_{n\in I_k}\lambda_n(T).
\]
With $M=C_0\norm T_{\mathcal W}$, by \eqref{eq:weak-eig},
\[
 |b_k|\le 2^{-k}\sum_{n\in I_k}\frac{M}{n}\le M2^{-k}.
\]
For the operator $B_b$ in \eqref{eq:dyadicspace}, write $\beta_n=b_k$ for $n\in I_k$. Then $|\beta_n|\le2M/n$, and the decreasing rearrangement of $(|\beta_n|)_{n\ge1}$ satisfies the same bound, as shown after \eqref{eq:dyadicspace}.

Let $2^k\le N<2^{k+1}$. Since $\sum_{n\in I_j}\beta_n=\sum_{n\in I_j}\lambda_n(T)$ for each $j$,
\begin{align}\label{eq:blockdifference}
 \left|\sum_{n=1}^N(\beta_n-\lambda_n(T))\right|
 &=\left|\sum_{n=2^k}^N(\beta_n-\lambda_n(T))\right|\notag\\
 &\le\sum_{n=2^k}^{2^{k+1}-1}(|b_k|+|\lambda_n(T)|)\le2M.
\end{align}
Let $(\widetilde\beta_n)_{n\ge1}$ be the eigenvalues of $B_b$ ordered by nonincreasing modulus. For either sequence $v=(\beta_n)$ or $v=(\widetilde\beta_n)$, $|v_n|\le2M/n$. For $N\ge2$,
\begin{align*}
 \left|\sum_{n=1}^Nv_n-\sum_{|v_n|>1/N}v_n\right|
 &\le\sum_{\substack{n\le N\\|v_n|\le1/N}}|v_n|
 +\sum_{\substack{n>N\\|v_n|>1/N}}|v_n|\\
 &\le1+2M\log\max\{1,2M\}.
\end{align*}
Indeed, the indices in the second sum satisfy $N<n<2MN$. The threshold sum is independent of the ordering, so
\begin{equation}\label{eq:reorder}
 \left|\sum_{n=1}^N\widetilde\beta_n-\sum_{n=1}^N\beta_n\right|
 \le2+4M\log\max\{1,2M\}.
\end{equation}
If $M=0$, both sequences vanish. Dividing \eqref{eq:blockdifference} and \eqref{eq:reorder} by $\log N$, we conclude that
\begin{equation}\label{eq:dyadic}
 \tau_X(T)=\tau_{\ell_2}(B_b)=\varphi\left(\left(2^{-k}\sum_{n\in I_k}\lambda_n(T)\right)_{k\ge0}\right),
\end{equation}
where $\varphi\colon Z\to\C$ is the linear functional $\varphi(b)=\tau_{\ell_2}(B_b)$.

Let $X,Y$ be real Banach spaces and $T\in\Lop(X,Y)$, with $T_{\C}$ defined by \eqref{eq:complex-op}. The norm on $X_{\C}$ is given by \eqref{eq:complex-norm}.
The canonical complex-linear map
\[
 J\colon(\ell_2(\R))_{\C}\to\ell_2(\C),\qquad J(x+iy)=(x_n+iy_n)_{n\ge1},
\]
is an isomorphism: by \eqref{eq:complex-norm},
\begin{align*}
 \tfrac1{\sqrt2}\norm{J(x+iy)}_2
 &\le \max\{\norm x_2,\norm y_2\}\le\norm{x+iy}_{(\ell_2(\R))_{\C}}\\
 &\le\left(\norm x_2^2+\norm y_2^2\right)^{1/2}=\norm{J(x+iy)}_2.
\end{align*}
Given $A\in\Lop(\ell_2(\C),X_{\C})$ with $\norm A\le1$, write $A|_{\ell_2(\R)}=A_1+iA_2$, where $A_j\in\Lop(\ell_2(\R),X)$ and $\norm{A_j}\le1$ for $j=1,2$. For $n\in\N$ and $\varepsilon>0$, choose $F_j\in\mathcal F(\ell_2(\R),Y)$ such that
\[
 \rank F_j<n,\qquad \norm{TA_j-F_j}\le a_n(TA_j)+\varepsilon\le x_n(T)+\varepsilon.
\]
For any $U\in\Lop(\ell_2(\R),Y)$ and $v\in\ell_2(\C)$,
\[
 \norm{U_{\C}J^{-1}v}_{Y_{\C}}=\sup_{\theta\in\R}\norm{U\operatorname{Re}(e^{i\theta}v)}_Y
 \le\norm U\norm v_{\ell_2(\C)}.
\]
Put
\[
 F=( (F_1)_{\C}+i(F_2)_{\C})J^{-1}\colon\ell_2(\C)\to Y_{\C}.
\]
Since $\rank F\le2n-2$ and $A=((A_1)_{\C}+i(A_2)_{\C})J^{-1}$,
\begin{align*}
 a_{2n-1}(T_{\C}A)&\le\norm{T_{\C}A-F}\\
 &\le\norm{(TA_1-F_1)_{\C}J^{-1}}+\norm{(TA_2-F_2)_{\C}J^{-1}}\\
 &\le2x_n(T)+2\varepsilon.
\end{align*}
Letting $\varepsilon\downarrow0$ and taking the supremum over $A$, we obtain
\begin{equation}\label{eq:complexification}
 x_{2n-1}(T_{\C})\le2x_n(T),\qquad \norm{T_{\C}}_{\mathcal W}\le4\norm T_{\mathcal W}.
\end{equation}
Indeed, for $m\in\N$ and $n=\lfloor(m+1)/2\rfloor$,
\[
 m x_m(T_{\C})\le2n x_{2n-1}(T_{\C})\le4n x_n(T)\le4\norm T_{\mathcal W}.
\] For every real Banach space $X$, define
\begin{equation}\label{eq:realtrace}
 \tau_X\colon\mathcal W(X,X)\to\R,\qquad \tau_X(T)=\operatorname{Re}\tau_{X_{\C}}(T_{\C}).
\end{equation}
For real $S,T\in\mathcal W(X,X)$ and $a,b\in\R$, \eqref{eq:homog}, \eqref{eq:additivity} and \eqref{eq:realtrace} yield
\[
 \tau_X(aS+bT)=\operatorname{Re}\tau_{X_{\C}}(aS_{\C}+bT_{\C})
 =a\tau_X(S)+b\tau_X(T).
\]
For $F\in\mathcal F(X,X)$, $F_{\C}\in\mathcal F(X_{\C},X_{\C})$ and $\tau_X(F)=\operatorname{Re}\tau_{X_{\C}}(F_{\C})=0$. By \eqref{eq:tracebound} and \eqref{eq:complexification},
\[
 |\tau_X(T)|\le4C_0\norm T_{\mathcal W}.
\]
For $A\in\mathcal W(X,Y)$ and $B\in\Lop(Y,X)$ over $\R$, by \eqref{eq:trace}, \eqref{eq:complexification} and \eqref{eq:realtrace},
\begin{align*}
 \tau_X(BA)&=\operatorname{Re}\tau_{X_{\C}}(B_{\C}A_{\C})\\
 &=\operatorname{Re}\tau_{Y_{\C}}(A_{\C}B_{\C})=\tau_Y(AB).
\end{align*}
For the real operator $D$ in \eqref{eq:diag}, the operator $JD_{\C}J^{-1}$ is diagonal on $\ell_2(\C)$ with entries $1/n$. Consequently $\lambda_n(D_{\C})=1/n$ for every $n\in\N$, and \eqref{eq:normalization} holds over $\R$.

For $b\in Z$, the eigenvalue sequences of $B_{\bar b}$ and $B_b$ are conjugate up to ties in modulus. By \eqref{eq:ties} and \eqref{eq:singulartrace}, $\varphi(\bar b)=\overline{\varphi(b)}$, so
\[
 \operatorname{Re}\varphi(b)=\tfrac12\big(\varphi(b)+\varphi(\bar b)\big)=\varphi(\operatorname{Re}b).
\]
Let $\varphi_{\R}\colon Z_{\R}\to\R$ be the restriction of $\varphi$. By \eqref{eq:dyadic} and \eqref{eq:realtrace}, for $T\in\mathcal W(X,X)$ over $\R$,
\[
 \tau_X(T)=\varphi_{\R}\left(\left(2^{-k}\sum_{n\in I_k}\operatorname{Re}\lambda_n(T_{\C})\right)_{k\ge0}\right).
\]
For $b\in Z_{\R}$, the diagonal operator in \eqref{eq:dyadicspace} over $\R$ satisfies $\varphi_{\R}(b)=\tau_{\ell_2(\R)}(B_b)$ by \eqref{eq:realtrace}.

\end{proof}

\section{Proof of Theorem B}
\begin{proof}
For complex $X$, the estimates \eqref{eq:weak-eig}, \eqref{eq:harmonic-sum} and \eqref{eq:ties} show that the sequence in \eqref{eq:quotient-trace} is bounded and that its class does not depend on the ordering at ties. For a bounded complex sequence $u=(u_N)_{N\ge1}$,
\begin{equation}\label{eq:ultrafilter-null}
 u\in c_0\quad\Longleftrightarrow\quad
 \lim_{N\to\mathcal U}u_N=0\quad\text{for every free ultrafilter }\mathcal U\text{ on }\N.
\end{equation}
Indeed, a null sequence has zero limit along every free ultrafilter. If $u\notin c_0$, there exists $\varepsilon>0$ such that $\{N:|u_N|\ge\varepsilon\}$ is infinite and belongs to some free ultrafilter $\mathcal U$. Then $|\lim_{N\to\mathcal U}u_N|\ge\varepsilon$.

For $S,T\in\mathcal W(X,X)$ and $\alpha,\beta\in\C$, let
\[
 u_N=\frac1{\log N}\sum_{n=1}^N
 \big(\lambda_n(\alpha S+\beta T)-\alpha\lambda_n(S)-\beta\lambda_n(T)\big)
 \qquad(N\ge2).
\]
Theorem~\ref{thm:weyltrace} was proved for an arbitrary free ultrafilter $\mathcal U$; in particular, the trace in \eqref{eq:singulartrace} is linear for each such $\mathcal U$. Hence $\lim_{N\to\mathcal U}u_N=0$ for every free ultrafilter. By \eqref{eq:ultrafilter-null}, $(u_N)_{N\ge2}\in c_0$, proving
\[
 \Psi_X(\alpha S+\beta T)=\alpha\Psi_X(S)+\beta\Psi_X(T).
\]
Likewise, if $A\in\mathcal W(X,Y)$ and $B\in\Lop(Y,X)$, applying \eqref{eq:ultrafilter-null} to the difference of the normalized eigenvalue sums for $BA$ and $AB$ proves $\Psi_X(BA)=\Psi_Y(AB)$.
The quotient norm of a bounded sequence modulo $c_0$ is the limit superior of the absolute values of its terms. Thus
\begin{align*}
 \norm{\Psi_X(T)}_{\ell_\infty(\C)/c_0(\C)}
 &=\limsup_{N\to\infty}\left|\frac1{\log N}\sum_{n=1}^N\lambda_n(T)\right|\\
 &\le C_0\norm T_{\mathcal W}\lim_{N\to\infty}\frac{1+\log N}{\log N}
 =C_0\norm T_{\mathcal W}.
\end{align*}
For a finite-rank operator the eigenvalue sum is eventually constant, so $\Psi_X(F)=0$. By \eqref{eq:normalization}, $\Psi_{\ell_2}(D)=[\mathbf 1]$.

For real $X$, the estimate \eqref{eq:complexification} and the complex case give
\[
 \Psi_X(T)=\operatorname{Re}\Psi_{X_{\C}}(T_{\C}),\qquad
 \norm{\Psi_X(T)}_{\ell_\infty(\R)/c_0(\R)}\le4C_0\norm T_{\mathcal W},
\]
where the real part of a complex sequence class is the class of the real parts of its terms. Complexification commutes with real linear combinations and compositions and preserves finite rank. The remaining identities follow from the complex case, including the normalization since $D_{\C}$ is similar to the complex diagonal operator with entries $1/n$.
\end{proof}

For any extended limit $\omega\colon\ell_\infty(\C)\to\C$, let $\bar\omega\colon\ell_\infty(\C)/c_0(\C)\to\C$ be the induced functional, $\bar\omega([u])=\omega(u)$. Since $\omega$ is positive and $\omega(\mathbf 1)=1$, it has norm one; so does $\bar\omega$. Define, for complex $X$,
\[
 \tau_X^\omega\colon\mathcal W(X,X)\to\C,\qquad
 \tau_X^\omega(T)=\bar\omega(\Psi_X(T)).
\]
For real $X$ use the restriction of $\bar\omega$ to $\ell_\infty(\R)/c_0(\R)$, which is real-valued. By Theorem~\ref{thm:quotient}, these maps are linear and satisfy the cyclic identity, vanish on finite-rank operators, and have the bound $|\tau_X^\omega(T)|\le C_{\K}\norm T_{\mathcal W}$ and the normalization $\tau_{\ell_2}^\omega(D)=1$. To see that these traces are spectral in the sense of \eqref{eq:spectral-def}, let $b=(b_k)_{k\ge0}$ be the dyadic means of the eigenvalues of $T$ as in the proof of Theorem~\ref{thm:weyltrace}. Equations \eqref{eq:blockdifference} and \eqref{eq:reorder} show that $\Psi_X(T)=\Psi_{\ell_2}(B_b)$ in the complex case, so the required linear functional is $b\mapsto\bar\omega(\Psi_{\ell_2}(B_b))$ on $Z$. For real $X$, take $b_k=2^{-k}\sum_{n\in I_k}\lambda_n(T_{\C})$ and the real diagonal operator $B_{\operatorname{Re}b}$. The block and rearrangement estimates apply to the real diagonal entries, since $|\operatorname{Re}b_k|\le |b_k|$. On the complex weak trace class these traces agree with the Dixmier traces of \cite{Ponge2023}.

\section*{Declarations}
\noindent\textbf{Funding.} G. Ara\'ujo is supported by Grant No. 2457/2026 from the Para\'iba State Research Foundation (FAPESq). G. Ara\'ujo, D. Pellegrino, A. Raposo Jr. and J. Santos are supported by CNPq (Brazil) through Grants No. 406457/2023-9 (CNPq/MCTI Call No. 10/2023) and No. 403964/2024-5 (MCTI/CNPq Call No. 16/2024). D. Pellegrino is also supported by CNPq Grant No. 305807/2025-0, and J. Santos by CNPq Grant No. 305655/2025-6.

\medskip
\noindent\textbf{Use of generative AI.} ChatGPT (OpenAI) was used to help develop the authors' initial outline of the proof, particularly in the choice of the Fourier coefficient and in the control of the harmonicity defect. It also assisted in an extensive search of the literature and in identifying supporting mathematical results. The authors take full responsibility for the mathematical content of the article.

\medskip
\noindent\textbf{Competing interests.} The authors declare that they have no competing interests.

\medskip
\noindent\textbf{Data availability.} No datasets were generated or analysed in this study.

% The authors should confirm the funding, competing interests and AI declarations before submission.

\end{document}